\documentclass[10pt]{article}
\usepackage[cp866]{inputenc} 

\begin{document} 

\centerline{\Large\bf Hidden unique possibilities  }
\centerline{\Large\bf of mathematical physics equations. } 
\centerline{\large\bf (Formalism of skew-symmetric forms)}
\centerline {\bf L.I. Petrova}
\centerline{\it Lomonosov Moscow State University, Faculty of Computational }
\centerline{\it Mathematics and Cybernetics, e-mail: ptr@cs.msu.su} 

\renewcommand{\abstractname}{Abstract}
\begin{abstract} 

It is shown that mathematical physics differential equations have properties that allow describing processes such as the structures emergence, discrete transitions, quantum jumps. The peculiarity is that such properties are hidden. They do not follow directly from the mathematical physics equations but are realized discretely in the solving process. This is due to the mathematical physics equations integrability, which, as shown, can be realized only discretely in the presence of any degrees of freedom. In this case, a transition occurs from the original coordinate space with a solution that is not a function (the solution derivatives do not compose a differential) to integrable structures with a solution that is a discrete function. It is the double solutions and spatial transitions that can describe the processes of the emergence of any structures or phenomena. 

Due to hidden properties, the mathematical physics equations have unique possibilities in describing physical processes and phenomena that cannot be described in the framework of other mathematical formalisms. 

Such results were obtained using skew-symmetric differential forms. 
\end{abstract}

{\bf Keywords}: the mathematical physics equations integrability; hidden properties; skew-symmetric differential forms; non-integrable manifolds; degenerate transformations.   

\section{Introduction} 

This paper deals with partial differential equations of mathematical physics describing continuous media (such as thermodynamic, gas-dynamic, cosmic), charged particles systems, and so on. It turns out that the equations of 
mathematical physics can not only describe evolutionary processes and changes in physical quantities (such as, for example, energy, pressure, density), but have the ability to describe discrete processes, such as the emergence of various structures, quantum jumps, the emergence of observable formations ( waves, vortices, turbulent pulsations) and so on.

However, such opportunities are hidden. They appear only in the process of equations solving. And this is due to the mathematical physics equations integrability, which is realized discretely. 

The mathematical physics equations, describing real processes, on the original tangent space turn out to be non-integrable. In this case, the equation solution is not a function: the derivatives of this solution do not compose a differential. 

The integrability of the equation can be realized only discretely in the process of solving the equations in the presence of any degrees of freedom. In this case, the resulting solution is a discrete function. 

Such features of mathematical physics equations integrability, as will be shown, reveal the unique possibilities of mathematical physics equations. 

\bigskip 
Section 2 notes the currently existing approaches to equations solving, which depend on the equations integrability, and the questions that arise in this case. 

At the same time, the characteristic properties of integrable equations are given, namely, the presence of differentials and relations in the differentials, which can be directly integrated. Such properties of integrable equations and the process of realizing integrability can only be described by skew-symmetric differential forms. 

Section 3 presents some properties and features of exterior and evolutionary skew-symmetric differential forms necessary for studying the differential equations integrability. 

Section 4 investigated the mathematical physics differential equations integrability, which depends on the described functions derivatives consistency and on the equations consistency if the equations of mathematical physics are a system of equations. (It should be noted that this is practically not taken into account when solving differential equations.) 

It is shown that the differential equations integrability is realized discretely. In this case, the hidden properties of mathematical physics equations are manifested, such as invariance, double solutions, discrete transitions, and so on.

It is these properties of the mathematical physics equations, as shown in Section 5, that reveal the unique possibilities of the mathematical physics equations in describing discrete processes and phenomena.

Section 5 summarizes the published works of the author, which describe such processes as the onset of vorticity and turbulence, discrete transitions, the implementation of Hamiltonian systems,
the connection of the equations of field theory with the equations of mathematical physics, and other results that turned out to be possible to obtain only thanks to the hidden properties of the mathematical physics equations.

\section {Existing approaches to solving differential equations of mathematical physics} 

We can mention analytical and numerical methods for solving of mathematical physics differential equations.

The use of analytical methods is possible if the differential equations can be reduced to an integrable form or if the integrability conditions are satisfied. 

Since the use of analytical methods even in simple cases is not always possible, many numerical methods were developed for the solution of mathematical physics equations.

The question arises, what is the correspondence between the solutions obtained by analytical and numerical methods and what is the accuracy of the studied real phenomena description by mathematical physics equations.

As shown in this work, there is a formalism of skew-symmetric differential forms, which allows not only to answer this question, but also to reveal the physical meaning of analytical and numerical solutions of mathematical physics equations.

\bigskip
Questions about solving of mathematical physics differential equations are related to the problem of integrability of differential equations. 

Obviously, the characteristic property of integrable equation is the presence of differentials and the ability to reduce the equation to the identical relation in the differentials, which will allow explicit integration of the differential equation.

Skew-symmetric differential forms  have such properties. Unlike other existing mathematical formalisms 
they can take the form of differentials and differential expressions. 

As you know, closed exterior forms are differentials. 
In addition, as established by the author, there are skew-symmetric differential forms that are obtained from differential equations and that generate closed exterior forms, i.e., differentials. Such skew-symmetric forms, which are evolutionary, made it possible to study the integrability and properties of mathematical physics differential equations.

This article is based on the properties of external and evolutionary skew-symmetric forms and on the features of their mathematical apparatus, which has non-traditional elements, such as non-identical relations, degenerate transformations, discrete transitions. Such a mathematical formalism is given only in the works of the author. Since this is not generally known, then
for the presentation of the results obtained with the help of this formalism, it is necessary to repeat some of the provisions of this formalism. This is presented in section 3. (For more details on skew-symmetric exterior and evolutionary forms, see [1].)
 
\section {Some properties and features of exterior and evolutionary skew-symmetric differential forms} 

Exterior differential forms (skew-symmetric forms on integrable manifolds and structures) have the properties, which underlie many mathematical formalisms. Historically, however, they are primarily used in differential geometry and topology. (It can be noted that the Frobenius theorem on the integrability conditions of the Pfaffian system is related to the skew-symmetric forms properties.) 

The exterior differential form of the degree $p$ ($p$-form) can be written
as [1, 2]
$$
\theta^p=\sum_{i_1\dots i_p}a_{i_1\dots i_p}dx^{i_1}\wedge
dx^{i_2}\wedge\dots \wedge dx^{i_p}\quad 0\leq p\leq n\eqno(1)
$$
Here $a_{i_1\dots i_p}$ is the function of independent variables $x^1, ..., x^n$,
$n$ is the space dimension, and
$dx^i$, $dx^{i}\wedge dx^{j}$, $dx^{i}\wedge dx^{j}\wedge dx^{k}$, \dots\
is the local basis subject to the skew-symmetry condition:
$$
\begin{array}{l}
dx^{i}\wedge dx^{i}=0\\
dx^{i}\wedge dx^{j}=-dx^{j}\wedge dx^{i}\quad i\ne j
\end{array}\eqno(2)
$$ 
The exterior form differential $\theta^p$ is expressed by the formula
$$
d\theta^p=\sum_{i_1\dots i_p}da_{i_1\dots
i_p}\wedge dx^{i_1}\wedge dx^{i_2}\dots \wedge dx^{i_p} \eqno(3)
$$
Since exterior forms are defined on integrable manifolds or structures, the exterior forms differential does not contain the basis differential from (since it is equal to zero). Such skew-symmetric forms can be closed.

\subsection {Closed exterior differential forms}

In mathematical formalisms the closed differential forms with
invariant properties appear to be of greatest practical utility. 

A form is 'closed', if its differential is equal to zero:
$$
d\theta^p=0\eqno(4)
$$
From condition (4) one can see that a closed form is a conservative 
quantity. (This means that it corresponds to a conservation law, 
namely, to some conservative physical quantity, i.e., quantity conservative under all invariant transformations.) 

If the exterior form is a closed inexact form, i.e., it is defined only on some structure (which is a pseudostructure by its metric properties), the closure condition is written as
$$
d_\pi\theta^p=0\eqno(5)
$$
And the pseudostructure $\pi$ obeys the condition
$$
d_\pi{}^*\theta^p=0\eqno(6)
$$
where ${}^*\theta^p$ is a dual form. I.e., the dual form describes a pseudostructure, on which a closed inexact exterior form is defined. In this case the interior differential rather than total one becomes equal to zero.

Closed inexact exterior form and the corresponding dual form compose a differential-geometric structure that describes conservative object, i.e., this quantity is conserved on a pseudostructure. (This can also correspond to some conservation law, i.e. to a conservative object.)

Following properties of closed exterior differential forms and features of their mathematical apparatus is related to the integrability of partial differential equations. 

1)	From the condition for the exterior form closedness 
$$
d\theta^p=0\eqno(7)
$$
it follows that the closed form is  {\it a differential}. 

The differential of the form is a closed form, since 
$$
dd\phi=0\eqno(8)
$$
where $\phi$ is an arbitrary exterior form. 

A closed exterior inexact form is an interior differential (on a pseudostructure). 
   
2) Since the closed form is a differential, this means that the closed form proves to be an invariant under all transforms that conserve the differential. 

The closed inexact exterior form and relevant dual form, composing the differential-geometrical structure, describe an invariant object.

3) The transformations of exterior skew-symmetric forms, which conserves the differential, are non-degenerate transformations. The unitary transforms (O-form), the tangent and canonical transforms (1-form), the gradient and gauge transforms (2-form) and so on are  examples of such transforms. 

4) Identical relations are an important element of the mathematical formalism of closed exterior forms. The identical relations express the fact that each closed exterior form is a differential of some exterior form. Such an identical relation can be written as
$$
d \phi=\theta^p\eqno(9)
$$
In this relation the form in the right-hand side has to be a {\it closed } one.   

It can be noted that practically in all branches of physics,
mechanics, thermodynamics, there are such identical relations. For example, the Poincare  invariant $ds\,=\,-H\,dt\,+\,p_j\,dq_j$, the second thermodynamics principle  $dS\,=\,(dE+p\,dV)/T$, and so on..  

It can be seen that closed exterior forms, firstly, are differentials, which indicates their integrability, and, secondly, the closed exterior forms consist of identical relations that can be integrated.

This indicates that if we understand the connection between differential equations and closed exterior forms, then we can solve the problem of differential equations integrability. It turns out that the evolutionary skew-symmetric forms properties reveal the connection of differential equations with closed exterior forms.

\subsection {Properties and features of evolutionary skew-symmetric differential forms.}

Evolutionary skew-symmetric forms [1], in contrast to exterior forms, are defined on nonintegrable manifolds (such as tangent manifolds of differential equations, Lagrangian manifolds and so on).
  
The evolutionary form can be written in a manner similar for exterior differential form [1]. However, in distinction from the exterior form differential, an additional term  will appear in the evolutionary form
differential. This is due to the fact that the evolutionary form basis
changes since such a form is defined on nonintegrable manifold.

The evolutionary form differential takes the form
$$
d\theta^p{=}\!\sum_{i_1\dots\i_p}\!da_{i_1\dots\i_p}\wedge
dx^{i_1}\wedge dx^{i_2}\dots \wedge
dx^{i_p}{+}\!\sum_{i_1\dots\i_p}\!a_{i_1\dots\i_p}
d(dx^{i_1}\wedge dx^{i_2}\dots \wedge dx^{i_p})\eqno(10)
$$
where the second term is connected with the basis differential being 
nonzero:
$d(dx^{i_1}\wedge dx^{i_2}\wedge \dots \wedge dx^{i_p})\neq 0$. 
(For the exterior form defined on integrable manifold one has
$d(dx^{i_1}\wedge dx^{i_2}\wedge \dots \wedge dx^{i_p})=0$). 

[Hereinafter a summing symbol $\sum$ and a symbol of exterior
multiplication $\wedge$ will be omitted. Summation over repeated indices
is implied.] 

\bigskip 
 
An evolutionary form, in contrast to exterior forms, cannot be a closed form, because its differential is not zero, since the second term of the evolutionary form differential associated with the basis differential is not zero.

The second term in the expression for the skew-symmetric 
form differential connected with the basis differential is expressed in
terms of the metric form commutator. For example, let us consider the first-degree form $\omega=a_\alpha
dx^\alpha$. 

The differential of this form can be written as
$d\omega=K_{\alpha\beta}dx^\alpha dx^\beta$ 
where $K_{\alpha\beta}$ is the evolutionary form  commutator. If we use the connectednesses  $\Gamma^\sigma_{\beta\alpha}$ to describe a nonintegrable manifold, then the commutator $K_{\alpha\beta}$ can be written in the form:
$$
K_{\alpha\beta}=\left(\frac{\partial a_\beta}{\partial
x^\alpha}-\frac{\partial a_\alpha}{\partial
x^\beta}\right)+(\Gamma^\sigma_{\beta\alpha}-
\Gamma^\sigma_{\alpha\beta})a_\sigma\eqno(11)
$$
The commutator $K_{\alpha\beta}$ is nonzero as the coefficients $a_{\mu }$   not potential, and connectednesses of 
a nonintegrable manifold is not symmetric [3].  
Thus it turns out that the commutator of the form $\omega$  and its differential does not equal to zero. 

Since the evolutionary form is not closed, then it cannot be a differential, like a closed exterior form. 
The evolutionary form has an unconventional mathematical apparatus containing non-degenerate self-varying relations and degenerate transformations. 

\subsection* {Non-identical self-varying relations.}

It was shown that identical relations lie at the
basis of exterior differential forms  the mathematical apparatus. 
In contrast to this, {\it nonidentical relations} lie at the basis of evolutionary differential forms 
mathematical apparatus.  Nonidentical relation appear in descriptions of any 
processes. 

Nonidentical relation may be written as well as the identical relation
$$
d\phi=\theta^p\eqno(12)
$$
However, on the right side of this relation there is an evolutionary skew-symmetric form, which is not a differential. 

The evolutionary relation turns out to be non-identical, since there is a differential on the left, and a skew-symmetric   on the right, which is not a differential.

An evolutionary nonidentical relation is self-varying, because, it is an evolutionary relation, and contains two members, one of which is immeasurable and cannot be compared with the other in the evolutionary process.
 
\subsection* {Degenerate transformation.}

One of the main elements of the mathematical formalism of the evolutionary forms theory is a degenerate transformation.  

Unlike a non-degenerate transformation of closed exterior forms, which conserves the differential, a degenerate transformation is a transformation that does not conserve the differential. (It can be noted that an example of a degenerate transformation is the Legendre transformation. However, it is practically applied only tacitly.) 

The degenerate transformation reveals the unique property of the evolutionary form. An evolutionary form can generate closed exterior forms that have invariant properties and correspond to the differential equations integrability.

As will be shown below, the evolutionary forms formalism, which has non-traditional elements, such as non-identical self-varying relations and degenerate transformations, makes it possible to discover the hidden properties of the mathematical physics equations, solve the problem of differential equations integrability and reveal their unique capabilities.

\section {On integrability of partial differential equations} 

As noted above, the characteristic property of integrable equations is the existence of differentials and identity relations in differentials, which allows us to explicitly integrate the differential equation. As has been shown, such property is  possessed by closed exterior forms. The differential equations integrability problem can be solved if one finds the connection of differential equations with closed exterior forms and understands how such a connection is realized. It is possible to solve this with the help of evolutionary skew-symmetric forms, which, as has been shown, can generate closed exterior forms.
\bigskip 

A qualitative study of differential equations has shown that the integrability of mathematical physics differential equations depends on the consistency of the derivatives of described functions  and on the consistency of equations if the equations of mathematical physics are the system of equations. 

\subsection {Investigation of the correspondence between derivatives in differential equations.} 

The correspondence between the derivatives of the sought-for functions in differential equations can be traced by the example of a first-order partial differential equation: 
$$ F(x^i,\,u,\,p_i)=0,\quad p_i\,=\,\partial u/\partial x^i \eqno(13)$$
Obviously, derivatives are consistent if they constitute a differential, i.e., the following relation is fulfilled 
$$ du\,=\,\theta\eqno(14)$$ 
where $\theta\,=\,p_i\,dx^i$ (summation over repeated indices is meant).  

In the general case, it does not follow (explicitly) from Eq. (13) that the derivatives $p_i\,=\,\partial u/\partial x^i $, satisfying the equation, constitute a differential. 

Derivatives $p_i$ can constitute a differential only if the conditions for the closedness of the form $\theta\,=\,p_jdx^j$ and the corresponding dual form are satisfied (in this case, the functional $F$ plays the role of dual form for $\theta $): 
$$\cases {dF(x^i,\,u,\,p_i)\,=\,0\cr
d(p_i\,dx^i)\,=\,0\cr}\eqno(15)$$

If we expand the differentials, then we get a system of homogeneous equations with respect to  $ dx ^ i $ and $ dp_i $ (in the $2n$ -dimensional cotangent space): 
$$\cases {\displaystyle \left ({{\partial F}\over {\partial x^i}}\,+\,
{{\partial F}\over {\partial u}}\,p_i\right )\,dx^i\,+\,
{{\partial F}\over {\partial p_i}}\,dp_i \,=\,0\cr
dp_i\,dx^i\,-\,dx^i\,dp_i\,=\,0\cr} \eqno(16)$$
The solvability condition of this system (i.e., 	vanishing of the determinant composed of coefficients at
$dx^i$, $dp_i$) has the form: 
$$
{{dx^i}\over {\partial F/\partial p_i}}\,=\,{{-dp_i}\over 
{\partial F/\partial x^i+p_i\partial F/\partial u}} \eqno (17)
$$ 
This condition determine an integrable structure, namely, a
pseudostructure (in its metric properties),  on which the form $\theta \,=\,p_i\,dx^i$ turns out to
be closed, i.e. it becomes a differential.  

It turns out that conditions (15) are satisfied only on integrable structures of the cotangent
space. In this case, on integrable structures (pseudostructures) derivatives of the differential equation constitute the differential $\delta u\,=\,p_idx^i\,=\,du$. That is, the solution becomes a function (discrete, since it is defined only on structures). On the cotangent space differential equation (13) turns out to be an integrable equation. 

But on the original the tangent space the differential equation (13) turns out to be non-integrable. In this case, the derivatives of the differential equation do not constitute the differential. Relation (14) turns out to be non-identical. On the right is a skew-symmetric form that is not a differential. Such a relation cannot be exactly integrated.
In this case, the solution of a differential equation is not a function. 

It turns out that equation (13) have double solutions, namely, on the original coordinate space and on the integrable structures. It can be seen that the transition from the first solution to the second is associated with the transition of the differential equation from the tangent space to the cotangent one. And this is due to the implementation of the integrability of the differential equation. Such a transition is described by a degenerate transformation.

\subsection {On the integrability of mathematical physics equations.} 

Equations of mathematical physics are practically a system of equations, i.e., they consist of several differential equations. Examples of mathematical physics equations are equations describing material media, such as thermodynamic, gas-dynamic, cosmic, and so on, which consist of the energy, momentum, angular momentum, and mass conservation laws equations. 

The system of equations integrability, an example of which are mathematical physics equations, depends on the consistency of the equations that make up the given system of equations.

The integrability mathematical physics equations was studied using the example of the system of equations consisting  of conservation laws equations. 

\subsection* {On the consistency of the conservation laws equations. Evolutionary relation.}

When investigating the consistency of the equations that form the system of equations, two coordinate systems were used: inertial and accompanying - associated with the manifold formed by the particles trajectories. (An example can be the Euler and Lagrangian coordinate systems.) 

In addition, the conservation laws equations were transformed into equations for functionals that characterize the described medium state. Functionals such as wave function, entropy, action functional, Poynting vector, Einstein tensor, etc. (which are also field theory functionals) are such state functionals [4]. 

When investigating the consistency of the conservation laws, we obtain an evolutionary relation in skew-symmetric forms for the state functional, which reveals the features of mathematical physics equations integrability. 

Let consider the consistency of the energy and momentum equations. In the accompanying coordinate system, the energy and momentum equations can be written in the form:
$$
{{\partial \psi }\over {\partial \xi ^1}}\,=\,A_1 \eqno(18)
$$ 
$$
{{\partial \psi}\over {\partial \xi^{\nu }}}\,=\,A_{\nu }\eqno(19)
$$

Here $\psi $ is the state functional,  $\xi^1$ and $\xi ^{\nu }$ are the coordinates along the trajectory and  in the direction normal to the trajectory, $A_1$ is a quantity that depends on specific features described material medium and
on external energy actions onto the medium, and $A_{\nu }$ are the quantities that depend on the specific
features of material medium and  external force actions.  

Here you can see the specific feature. For the same quantity (functional $\psi$), there are two equations. How can such a system be investigated? 

Equations (18) and (19) can be folded into a relation
$$
d\psi\,=\,A_{\mu }\,d\xi ^{\mu }\eqno(20)
$$ 
Relation (19) can be written as
$$
d\psi \,=\,\omega\eqno(21) 
$$
where $\omega \,=\,A_{\mu }\,d\xi ^{\mu }$ is the skew-symmetric
differential form of the first degree. 

Since the conservation laws equations are evolutionary, the resulting relation and the skew-symmetric form $\omega$ will also be evolutionary. 

If we also take into account the angular momentum and mass conservation laws, then the evolutionary relation takes the form
$$
d\psi \,=\,\omega^p\eqno(22)
$$ 
where the degree of the form  $p$ takes the values $p\,=\,0,1,2,3$.  

(The article [5] describes the derivation of the evolutionary relation for a gas-dynamic medium. This evolutionary relation is given in section 5.) 

\subsection*{Properties and features of the evolutionary skew-symmetric form and evolutionary relation.}

An evolutionary skew-symmetric form, unlike exterior forms, cannot be closed. 

As shown, the evolutionary relation (21) and evolutionary form $\omega$ were obtained in the accompanying coordinate system associated with the manifold formed by the described medium particles trajectories. The manifold formed by the particles trajectories is deformable, i.e., a non-integrable manifold. As shown in subsection 3.2 evolutionary form defined on the non-integrable manifold cannot be closed because its commutator (which contains a non-zero metric form commutator), and, consequently, the differential is not equal to zero.
  
Since the evolutionary form $\omega$ is not a differential, the evolutionary  relation turns out to be {\it non-identical (the differential is on the left, and the skew-symmetric form, which is not a differential, is on the right)}. 

\bigskip 

In a similar way, it can be shown that relation (22) will also be nonidentical. In this case, the evolutionary forms of degrees $ 1-3 $ will be non-closed, since the commutator of these evolutionary forms will contain nonzero commutators of the metric form of the first, second, and third degrees, which determine the torsion, rotation, and curvature, respectively.(It should be noted that the resulting evolutionary relation remains nonidentical regardless of the accuracy with which were written the conservation laws equations. And it has a physical meaning.) 
\bigskip

Nonidentical evolutionary relation turns to be self-changing relation as it evolutionary and comprises members, which do not agree between themselves.  

\subsection* {Non-integrability of mathematical physics equations on the original tangent space.} 

The nonidentity of the evolutionary relation points out to the fact that the 
conservation law equations appear to be inconsistent. This means that
the initial set of equations of mathematical physics proves to be nonintegrable.  

The nonidentical evolutionary relation cannot be explicitly integrated, since there is a not closed form on the right, which is not a differential. Mathematical physics equations cannot be folded into the identical relation and explicitly integrated. 
It turns out that on the original tangent space the mathematical physics equations turn out to be non-integrable. 

In this case, the solutions of mathematical physics equations will depend not only on variables. They will depend on the evolutionary form commutator. The solutions will not be analytical solutions, i.e., functions. The derivatives of such solutions do not form a differential. (As will be shown below, such solutions have a physical meaning.)

\subsection* {Realization of the mathematical physics equations integrability.}

The mathematical physics equations integrability can be realized only if the identical relation (which can be integrated) is obtained from the non-identical relation. 
To obtain the identical relation from an evolutionary non-identical relation, it is necessary to obtain the closed exterior form from the evolutionary form included in the evolutionary relation.

Since the differential of the evolutionary form is non-zero, and the differential of the closed exterior form is equal to zero, the transition from the evolutionary form to the closed exterior form is possible only with a degenerate transformation, i.e., with a transformation that does not conserve the differential.

\bigskip
Degenerate transformations can be done under additional conditions, which are resulted from the implementation of any degrees of freedom. Vanishing of functional expressions such as 
determinants, Jacobians, Poisson brackets, residues, and others correspond to these additional conditions. Vanishing of these functional expressions is a condition for the closedness of the dual form describing the implemented integrable structure.

Additional conditions can be realized upon self-change of a non-identical evolutionary relation. This can only happen discretely. 

The conditions for the degenerate transformation define integrable structures (pseudostructures): a characteristic (the determinant vanishes), singular points (the Jacobian is zero), etc. 

The degenerate transformation is realized as a transition from the original tangent, non-integrable  space to  a cotangent, integrable manifold. 

An example of the degenerate transformation is the Legendre transformation. However, the Legendre transform is usually performed implicitly and does not take into account the discrete nonequivalent transition from one space to another. (An example of such a transition is the transition from a Lagrangian manifold to a Hamiltonian one.)

\subsection*{Realization of a closed inexact exterior form. Obtaining an identical relation from a nonidentical.} 

When the condition of a degenerate transformation is fulfilled, the following transition is realized 

$d\omega^p\ne 0 \to $ (degenerate transformation) $\to d_\pi \omega^p=0$,
$d_\pi{}^*\omega^p=0$

Condition  $d_\pi{}^*\omega^p=0$ defines a closed dual form that describes some integrable structure $\pi$ (a pseudostructure in its metric properties). The condition $d_\pi \omega^p=0$  defines the closed inexact (on an integrable structure) exterior form. 

(The following should be emphasized here.
The realization of a closed external form, which is a differential and, therefore, an invarant, indicates that the mathematical physics equations have invariant properties.)

On an integrable structure, from  evolutionary relation (22) implies relation  
$$
d\psi_\pi=\omega^p_\pi\eqno(23)
$$ 
which turns out to be identical since the form $\omega^p_\pi$ is a differential (interior). 

Since the identity relation can be integrated, this indicates that the equation of mathematical physics turns out to be integrable on the integrable structure. I.e., the mathematical physics equation solution is a function. But such an analytical solution is a discrete function, since it is determined only on structure. This is the so-called generalized solution.

\bigskip 

So, it turns out that mathematical physics equations have double solutions, namely, on the original coordinate space and on the integrable structure. Moreover, the solution on the original coordinate space is not a function. And the solution on the integrable structure is a discrete function.

As shown, there is a connection between these solutions, which is implemented discretely. When the conditions for the degenerate transformation are satisfied, integrable structures are realized and a transition occurs from the original tangent space with a solution that is not a function to the integrable structures of the cotangent manifold with a generalized solution, which is a discrete function.

\subsection*{Physical meaning of double solutions.} 

It was shown that mathematical physics equations have double solutions: on a non-integrable initial coordinate space and on integrable structures. It makes physical sense. Inexact solutions describe the nonequilibrium state of the medium under study. Exact solutions (which are discrete functions) describe the local-equilibrium state of the medium.

This follows from the evolutionary relation $d\psi \,=\,\omega^p$ and is connected with the fact that the evolutionary relation contains the functional $\psi$, which characterizes the state of the medium. 

The presence of the state functional differential indicates the presence of a state function, which corresponds to the equilibrium state of the medium. But there is a subtlety here. Since the evolutionary relation is not identical, it is impossible to obtain the differential of the state functional $d\psi$ from it. This indicates the absence of a function of state and means that the state of the described medium is nonequilibrium. Inexact solutions describe such a nonequilibrium state of the medium. Obviously, the internal forces that cause non-equilibrium should be described by a commutator of the open evolutionary form $\omega^p$. 

From the identical relation $d\psi_\pi=\omega^p_\pi$, which is realized from the evolutionary relation, one can obtain the differential  $d\psi_\pi$ (interior, on an integrable structure). And this indicates the presence of a state function and the transition of the medium to a locally equilibrium state. 
(In this case, the general state of the medium remains non-equilibrium.) The locally-equilibrium state of the medium is described by a discrete function. 

The transition of medium from the nonequilibrium state to the local-equilibrium state means that an immeasurable quantity, which was described by the evolutionary form commutator and acted as an internal force, becomes a measurable quantity (potential, corresponding to the given medium). This manifests itself as the emergence of any discrete structures and observable formations, an example of which are waves, vortices, turbulent pulsations.
It should be emphasized once again that mathematical physics equations solutions on integrable structures are related to solutions on the original coordinate space. This must be taken into account when describing real physical processes by mathematical physics equations. 

(It should be noted that inconsistency of the conservation law equations reflects the non-commutativity of conservation laws. It is with the non-commutativity of conservation laws that the processes occurring in material environments are associated.)

\bigskip 

So, it turns out that mathematical physics equations have double solutions: on a non-integrable initial coordinate space and on integrable structures. 

How can you get these solutions?

Since mathematical physics equations on the original coordinate space have solutions that are not functions, such solutions can only be obtained by numerical methods.

And on integrable structures, solutions that are functions can be obtained using analytical methods. Such solutions can also be obtained by numerical methods. But in this case, the equations are modeled on the integrable manifold, and not on the original space. I.e., other, nonequivalent, coordinate systems are used. (It is seen that double solutions cannot be obtained using a unified continuous simulation of mathematical physics equations.)

As is known, at present there are both analytical and numerical methods for solving these equations. In this case, one and the same phenomenon is sometimes described both by numerical methods on the original coordinate space, and by analytical methods on integrable spaces, such as phase (covariant, cotangent) spaces (which are obtained from the integrability conditions imposed on the equations).

But, in this case, the solutions of mathematical physics equations describing the same phenomenon but obtained by different methods (numerical methods on the original coordinate space or analytical methods on integral structures), will be different. This does not mean that the solutions are wrong. The solutions of mathematical physics equations, obtained by both numerical and analytical methods, are correct, but in each case they describe only one of the sides of the physical phenomenon. Therefore, in order to obtain a complete description of a physical phenomenon, it is necessary to simultaneously use both methods for solving equations of mathematical physics and take into account the relationship between these approaches.

(It can be noted that in the continuum mechanics, describing the change in quantities characterizing material media, numerical methods are usually used in the initial coordinate space, and in the continuum physics, describing conserved (under non-degenerate transformations) physical quantities or objects, usually mathematical physics equations are solved by analytical methods or numerical methods, but on integrable manifolds. This can be explained by the fact that in continuum mechanics they are interested in evolutionary processes described by numerical solutions, and in continuum physics they are interested in conserved objects described by analytical solutions.)

\bigskip 

Investigation of the integrability of differential equations describing real processes and phenomena showed that the differential equations integrability is realized discretely in the presence of any degrees of freedom. At the same time, unique hidden properties of mathematical physics equations are manifested, such as invariance, double solutions, discrete transitions, and so on, which makes it possible to describe evolutionary processes, the emergence of physical structures, observed formations, and so on.
It should be emphasized once again that such results were obtained using skew-symmetric evolutionary forms, which are defined on non-integrable manifolds and have an unconventional mathematical apparatus containing nonidentical relations and degenerate transformations. 

Such properties of differential equations reveal the unique capabilities of mathematical physics equations in describing physical processes and phenomena that cannot be described within the framework of other mathematical formalisms.

\section {Unique hidden possibilities of the mathematical physics equations} 

In the author`s published works, which are briefly presented below, the properties and possibilities of some mathematical physics equations are given. The features of the described phenomena are shown, which are due to the hidden properties of the mathematical physics equations. 
 
\subsection* {Double Solutions of the Euler and Navier-Stokes Equations. Process of Origination the Vorticity and Turbulence} 

The Euler and Navier-Stokes equations (which describe a gas-dynamic medium, i.e., the ideal and viscous gas flow) are examples of mathematical physics equations. 

As is known, the proof of the existence and smoothness of the Navier-Stokes equations solution was declared as a millennium problem. 

But as follows from this article, and also as shown in [5], there is no smooth analytical solution. The Euler and Navier-Stokes equations, as well as the mathematical physics equations, have double solutions, namely, the solution on the original coordinate space and solution on integrable structures. In this case, the transition from the first solution to the second describes the transition of a gas-dynamic medium from a nonequilibrium state to a local-equilibrium state, which is accompanied by the emergence of vorticity (in an ideal gas) or the emergence of turbulence (in a viscous gas). 

This follows from the evolutionary relation, which is obtained from the Euler and Navier-Stokes equations. The properties of the Euler and Navier-Stokes equations can be found in [5].

It can be noted that the evolutionary relation obtained from the Euler and Navier-Stokes equations can appear as follows [5]: 
$$ds\,=\,A_{\mu} d\xi ^{\mu}$$ 
Here the coefficients $A_{1}$ obtained from the energy conservation law equation, respectively, are equal to 
$A_{1}=0$ for the ideal gas, and  
$$A_1\,=\,{1\over {\rho }}{{\partial }\over {\partial x_i}}
\left (-{{q_i}\over T}\right )\,-\,{{q_i}\over {\rho T}}\,{{\partial T}\over {\partial x_i}}
\,+{{\tau _{ki}}\over {\rho }}\,{{\partial u_i}\over {\partial x_k}}$$ 
for the viscous gas. And the coefficients $A_{\nu1}$ obtained from the momentum conservation law equation are [6]:  
$$A_{\nu}={{\partial h_{0}}\over {\partial \xi^{\nu}}}+{({u_1}^2+{u_2}^2)}^{1/2}\zeta-
F_{\nu}+{{\partial U_{\nu}}\over {\partial t}}$$      

It should be emphasized here that the entropy $s$, which is included in the evolutionary relation, depends on the space-time coordinates  $\xi^{\mu}$, in contrast to the thermodynamic entropy, which depends on thermodynamic variables. It is this entropy (and not thermodynamic) that describes the gas-dynamic medium state.

\subsection* {Connection of the field theory equations with the mathematical physics equations.} 

At the end of the 19th and the beginning of the 20th centuries, a situation arose that Henri Poincare called the "crisis of physics". A lot of phenomena and structures have accumulated that turned out to be impossible to describe by the mathematical physics equations. Therefore, in different areas of physics, new equations- the field theory equations were obtained (such as the equations of Einstein, Maxwell, Schr?dinger, Dirac, Heisenberg). These equations are based on the invariance and covariance properties, which are necessary to describe measurable physical structures and observable phenomena, and which, it was assumed, do not have the mathematical physics equations. 

It turns out that the equations of mathematical physics have invariant properties. But, as shown, such properties are hidden. They are realized discretely and this is to the realization of the mathematical physics equations integrability, which is described by the realization of closed skew-symmetric exterior and dual forms that have invariant and covariant properties. Invariant and covariant properties of mathematical physics equations reveal a connection of the field theory equations with the mathematical physics equations [7]. 

The evolutionary relation that is obtained from the mathematical physics equations, has the properties of the field theory equations. 

The work [7] shows the correspondence between the evolutionary relation and the field theory equations, which indicates the connection of the field theory equations with the of mathematical physics equations. 

As already been mentioned, the evolutionary relation is a relation for functionals such as the action functional, entropy, wave function, Lagrangian, Einstein tensor, Poynting vector, etc., which are also functionals of field theory [4]. 

The field theory equations also have the form of relations written in terms of skew-symmetric forms or their tensor or differential analogs. So 

-	the Einstein equation is a relation in skew-symmetric forms (or a tensor relation);

-	the Maxwell equations have the form of tensor relations; 

-	the Schrodinger equations have the form of relations expressed in terms of derivatives and their analogs. 

From the evolutionary relation closed inexact exterior forms can be obtained, which are solutions to field theory equations. 

Evolutionary relation of the mathematical physics equations unifies the field theory equations, reveals their internal connection, reveals properties common to all field theory equations. This can serve as an approach to a unified and general field theory.

\subsection* {Duality of conservation laws.} 

The mathematical physics equations reveal the conservation laws duality [8].

The hidden properties of the mathematical physics equations make it possible to understand the features of conservation laws, which have different meanings in thermodynamics, physics and mechanics. 

In mechanics and continuous media physics the concept of ``conservation laws" relates to the conservation laws for energy, linear momentum, angular momentum. These are conservation laws for material media. They are described by differential equations. 

In areas of physics related to the field theory and in the theoretical mechanics ``the conservation laws" is conservation laws, according to which the conservative physical quantities or objects exist. Such conservation laws are described by closed exterior and dual forms, which, as noted, are conserved quantities or objects. Conservation laws for physical fields are such conservation laws. 

The first principle of thermodynamics is an example of evolutionary relation. It is obtained from two conservation laws, namely, from the energy conservation law and the linear momentum conservation law.
The second principle of thermodynamics, which determines the thermodynamic entropy, is an example of the identity relation that is obtained from the evolutionary relation (the first principle of thermodynamics) when implementing the integrating factor: 1 / T, where T is the temperature.

The realization of closed exterior forms from the mathematical physics equations, which consist of equations for energy, linear momentum, angular momentum, reveals the connection between the conservation laws for physical fields and the conservation laws for material media. This connection is realized discretely in the evolutionary process [8]. Realization of closed exterior forms describing the conservation laws for physical fields from the mathematical physics equations describing material media, indicates that physical fields are generated by material media. 

\bigskip 

It should be emphasized here that there is no exact law of conservation of energy. There is no law of conservation of energy as a conserved quantity, which is described by a closed external form and corresponds to the conservation law for physical fields. But there is a balanced (differential) law of conservation of energy, which is described by the differential equation [8].

\subsection*{Hamiltonian systems realization process.} 

The mathematical physics equations reveal the integrable Hamiltonian systems realization process. 

The hidden properties of the mathematical physics equations make it possible to describe the process of Hamiltonian systems realization from differential equations. This is described in the article [9]. However, the article [9] describes not the properties of the Hamiltonian system and the Hamiltonian function (to which many papers are devoted), but shows how Hamiltonian systems are realized from differential equations, and describes the implementation process.

In classical mechanics, the Hamiltonian system is realized from the Euler-Lagrange equations as the Euler-Lagrange equation integrability [9]. 
It should be noted here that the Hamiltonian system can be obtained from the Euler-Lagrange equation if no additional conditions are imposed on it. In this case, the Euler-Lagrange equation turns out to be non-integrable on the tangent space. However, when some degrees of freedom are realized (when the constraints locally become holonomic), the integrability condition (one of the Hamiltonian system relations) is satisfied, and the Hamiltonian system for the Hamiltonian function is realized from the Euler-Lagrange equation. 

The Hamiltonian system is realized discretely. In this case, the transition from a tangent manifold with the Euler-Lagrange equation to sections of a cotangent manifold with a Hamiltonian system realized by the Legendre transformation, which is a degenerate  transformation (a transformation that does not preserve the differential).
\subsection*{The emergence of the physical structures and the observed formations.} 
The mathematical physics equations reveal the process of the emergence of physical structures and the emergence of observed formations [10]. 

At a degenerate transformation from the evolutionary form, the closed inexact exterior form and relevant dual form are realized, which form the differential geometrical structure, namely, the pseudostructure with conserved quantity. The differential-geometrical structure describes the physical structure, i.e., the structure on which the exact conservation law is fulfilled. The physical fields and relevant manifolds are formed by such physical structures [10]. Massless particles, structures made up by eikonal surfaces and wave fronts, and so on are examples of physical structures. 

The emergence of the physical structures is accompanied by the emergence of observable formations such as waves, vortices, fluctuations, turbulent pulsations, and so on in the described media [10]. 

The description of process the physical structures emergence  has features that are not present in any of the existing mathematical formalisms. The physical structures emergence is described by the transition from a tangent, non-integrable, space to a cotangent, integrable manifold (i.e., from one spatial object to another not equivalent spatial object). This can only be described by skew-symmetric forms using the degenerate transform. 

{\footnotesize It can be noted here that there is a connection between a degenerate transformation and a non-degenerate one. This can be shown on the Lagrangian formalism example. The degenerate transformation is a transition 
from coordinate space ($q_j,\,\dot q_j)$) to integrable structures cotangent manifold ($q_j,\,p_j$). And the nondegenerate transformation is a transition in cotangent manifold from some structure  ($q_j,\,p_j$) to another structure ($Q_j,\,P_j$). 

The degenerate transformations describe the process of emergence of physical structures, whereas the nondegenerate transformations describe transitions from one structure to another.}

\subsection*{Discrete transitions and quantum leaps. } 
The mathematical physics equations can describe discrete quantum transitions [11].

It should be emphasized that discrete processes can only describe the mathematical physics equations on which the integrability conditions are not imposed. It is the realization of the integrability conditions, which occurs discretely in the equations solving process, that reveals the hidden properties of the mathematical physics equations, which make it possible to describe discrete transitions, the emergence of various structures and observed formations. 

This also reveals the mechanism of occurrence of quantum leaps. In transition from one structure to another, the conserved quantity described by a closed exterior form and the pseudo structure described by the dual form, changes discretely. The discrete change in the conserved quantity and pseudo structure produces a quantum which is emerged in transition from one structure to another. Evolutionary form commutator, formed at the time of the structure emergence, determines this quantum characteristics. 

It can be emphasized once again that discrete transitions are described by the degenerate transformation, in which there is a transition from the original coordinate space to integrable structures.

\section{Conclusion} 

The unique hidden possibilities of the mathematical physics differential equations were discovered when studying the mathematical physics equations integrability. In this case, it was taken into account that the condition for integrability is the consistency of the desired function derivatives (i.e., the presence of a differential) and the consistency of equations if the mathematical physics equations are the system of equations. (These conditions are usually not taken into account in solving the mathematical physics equations.) 

It was found that the mathematical physics equations integrability is realized discretely in the equations solving process (and is not a constant property of the mathematical physics equations). At the same time, hidden properties of mathematical physics equations are revealed, which turn out to be unconventional. Firstly, there is a transition from the tangent, non-integrable space, on which the mathematical physics equations are defined, to the integrable structures of the cotangent manifolds . Secondly, on integrable structures, the solution of mathematical physics equations becomes a discrete function. This is the so-called generalized solution. In this case, the solution in the original coordinate space remains non-integrable, i.e., it is not a function. It turns out that the mathematical physics equations have double solutions. 

Such properties of the mathematical physics equations, as shown, allow us to describe discrete transitions, quantum leaps, the emergence of various structures and observed formations, such as waves, vortices, turbulent pulsations, and so on. Due to such possibilities, the mathematical physics equations reveal the mechanism of many ongoing processes and phenomena. This is demonstrated in the published works of the author, which are briefly presented in this paper. 

Such results were obtained by skew-symmetric differential forms, which are differentials and differential expressions, which allows solving integrability problems. At the same time, in addition to closed exterior forms that have invariant properties, the skew-symmetric forms obtained by the author were used, which are defined on non-integrable manifolds and have evolutionary properties. The evolutionary forms theory contains elements such as non-identical relations and degenerate transformations that are practically not explicitly applied in the existing mathematical formalisms.

[1] Petrova L.I. Role of skew-symmetric differential forms
in mathematics, 2010,  http://arxiv.org/pdf/1007.4757.pdf  

[2] Bott R., Tu L.~W., {\it Differential Forms in Algebraic Topology}, Springer, NY, 1982.  

[3] Tonnelat M.-A., {\it Les principles de la theorie electromagnetique
et la relativite}, Masson, Paris, 1959. 

[4]  Petrova L. Connection between functionals of the field-theory equations and state functionals 
of the mathematical physics equations. {\it Journal of Physics: Conference Series}, IOP Publishing ([Bristol, UK], England), 2018, Vol 1051, No. 012025, 1- 8. 

[5]  Petrova L.I. Double Solutions of the Euler and Navier-Stokes Equations. Process of Origination the Vorticity and Turbulences. {\it Fluid Mechanics}, Science Publishing Group, 2017, Mar. 21, 6-12. 

[6] Clark J F., Machesney M. {\it The Dynamics of Real Gases}, Butterworths, London, 1964. 

[7] Petrova L. Evolutionary Relation of Mathematical Physics Equations Evolutionary Relation as Foundation of Field Theory. Interpretation of the Einstein Equation, {\it Axioms}, 2021, MDPI (Basel, Switzerland), Vol 10(2), No.46, 1-10; https://doi.org/10.3390/axioms10020046
 
[8] Petrova L.I. Duality of conservation laws and their role in the processes of emergence of physical structures and formations, {\it Mathematics for Applications}, Brno University of Technology (Czech Republic), 2021, Vol 8, No.9, 55-70. 

[9] Petrova L, Qualitative Investigation of Hamiltonian Systems by Application of Skew-Symmetric Differential Forms, 
{\it Symmetry},  MDPI (Basel, Switzerland), 2021, Vol 13,  No.1, 1-7 

[10] Petrova L.I. Spontaneous emergence of physical structures and observable formations: fluctuations, waves, turbulent pulsations and so on, {\it Journal of Applied Mathematics and Physics},  Scientific Research Publishing (United States),  2016, Vol 4, No 5, 864-870. 

[11] Petrova L. Discrete Quantum Transitions, Duality: Emergence of Physical Structures and Occurrence of Observed Formations (Hidden Properties of Mathematical Physics Equations), {\it Journal of Applied Mathematics and Physics (JAMP)}, Scientific Research Publishing, 2020, Vol.8, No.9, 1911-1921. 

\end{document}